\documentclass[twocolumn]{autart}    

\usepackage{graphicx}          
\usepackage{lineno,hyperref}
\usepackage{amsmath}

\usepackage{amssymb}
\usepackage{amscd,multirow}
\usepackage{graphicx}
\usepackage{epstopdf}
\usepackage[titletoc]{appendix}
\usepackage{subfigure}
\usepackage[justification=centering]{caption}
\usepackage{cite,natbib}
\usepackage{wrapfig}
\usepackage{graphicx}
\usepackage{textcomp}
\usepackage{amsfonts}
\usepackage{subfigure}
\usepackage{subeqnarray,cases}
\usepackage{hyperref}
\usepackage[misc]{ifsym}
\usepackage{algorithm}
\usepackage{setspace}

\hypersetup{hypertex=true,
            colorlinks=true,
            linkcolor=blue,
            anchorcolor=blue,
            citecolor=blue}

\allowdisplaybreaks[4]

\newtheorem{theorem}{Theorem}
\newtheorem{corollary}{Corollary}
\newtheorem{lemma}{Lemma}
\newtheorem{proposition}{Proposition}
\newtheorem{remark}{Remark}

\begin{document}

\begin{frontmatter}

\title{Fast  primal-dual algorithm  via  dynamical system for a linearly  constrained convex optimization problem\thanksref{footnoteinfo}} 

\thanks[footnoteinfo]{ Corresponding author: Ya-Ping Fang.}

\author[Scu]{Xin He}\ead{hexinuser@163.com},    
\author[CUIT]{Rong Hu}\ead{ronghumath@aliyun.com},               
\author[Scu]{Ya-Ping Fang}\ead{ypfang@scu.edu.cn}  

\address[Scu]{Department of Mathematics, Sichuan University, Chengdu, Sichuan, P.R. China}  
\address[CUIT]{Department of Applied Mathematics, Chengdu University of Information Technology, Chengdu, Sichuan, P.R. China}             

\begin{abstract}                          
By time discretization of a second-order primal-dual  dynamical system with damping $\alpha/t$ where an inertial construction in the sense of Nesterov is needed only for the primal variable, we propose a fast  primal-dual algorithm for a linear equality constrained convex optimization problem.  Under a suitable scaling condition, we show that the proposed algorithm enjoys a fast convergence rate for the objective residual and the feasibility violation, and the decay rate can reach $\mathcal{O}(1/k^{\alpha-1})$ at the most.  We also study convergence properties of the corresponding primal-dual dynamical system to better understand the acceleration scheme. Finally, we  report numerical experiments to demonstrate the effectiveness of the proposed algorithm.
\end{abstract}

\end{frontmatter}

\section{Introduction}
Let $\mathbb{R}^{n}$ be an $n$-dimensional Euclidean space   with the  scalar product $\langle \cdot,\cdot\rangle$ and the corresponding induced norm $\|\cdot\|$.
 Let  $f: \mathbb{R}^{n}\to\mathbb{R}\cup\{+\infty\}$ be a  proper, lower semicontinuous and convex  function, $A\in\mathbb{R}^{m\times n}$,  and $b\in\mathbb{R}^{m}$.  Consider the linearly  constrained convex optimization problem
\begin{equation}\label{ques}
				\min_{x\in\mathbb{R}^{n}}  \quad f(x), \quad s.t.  \  Ax = b.
\end{equation}
The problem \eqref{ques} is a basic model for many important applications arising in various areas, such as compressive sensing, image processing, global consensus  and machine learning problems. See e.g. \cite{Candes2008,Boyd2011,Lin2020,Feijer2010,Zhang2010,Zhu2018IEEE,Wang2021}.

Denote the KKT point set of the problem \eqref{ques} by $\Omega$. For any  $(x^*, \lambda^*)\in\Omega$, we have
 \begin{equation}\label{saddle_point}
	\begin{cases}
		-A^T\lambda^* \in  \partial f(x^*),\\
		Ax^*=b,
	\end{cases}
\end{equation}
where
\[\partial f(x) = \{v\in\mathbb{R}^n\ |\ f(y)\geq f(x)+\langle v,y-x\rangle,\quad \forall y\in\mathbb{R}^n\}. \]
Recall the Lagrangian function of the problem \eqref{ques},
\begin{equation*}\label{AugL}
	\mathcal{L}(x,\lambda)= f(x)+\langle \lambda,Ax-b\rangle,
\end{equation*}
where $\lambda\in \mathbb{R}^{m}$ is the Lagrange multiplier. From \eqref{saddle_point} we have
\begin{equation*}\label{eq_Larg}
	\mathcal{L}(x^*,\lambda)\leq  \mathcal{L}(x^*,\lambda^*)\leq  \mathcal{L}(x,\lambda^*),\ \ \forall  (x,\lambda)\in \mathbb{R}^{n}\times\mathbb{R}^m.
\end{equation*}
Throughout this paper, we always assume  $\Omega\neq \emptyset$.
\subsection{Literature review}
A  benchmark algorithm for the problem \eqref{ques} is the augmented Lagrangian method (ALM)
\begin{eqnarray*}\label{al:in1}
	\begin{cases}
		x_{k+1}\in\mathop{\arg\min}_x f(x)+\langle \lambda_k,Ax-b\rangle+\frac{\sigma}{2}\|Ax-b\|^2,\\
		\lambda_{k+1} = \lambda_k+\sigma(Ax_{k+1}-b).
	\end{cases}
\end{eqnarray*}
ALM  plays a fundamentally theoretical and algorithmic role in solving the problem \eqref{ques}.  Here, we mention some of  nice works concerning fast convergence properties of  ALM and its variants. By applying Nesterov's acceleration  technique  \cite{Nesterov2018,Nesterov1983} to ALM, \cite{HeY2010} developed an accelerated augmented Lagrangian method (AALM) for the problem \eqref{ques} and  proved that AALM enjoys the  $\mathcal{O}(1/k^2)$ convergence rate when  $f$ is differentiable. When  $f$ is nondifferentiable,  the   $\mathcal{O}(1/k^2)$ convergence rate of AALM was established in \cite{Kang2013}. \cite{Kang2015} further proposed an inexact version of AALM and demonstrated the $\mathcal{O}(1/k^2)$ convergence rate   under the strong convexity assumption of $f$.  \cite{Huang2013} considered an accelerated linearized Bregman method for  solving the basis pursuit and related sparse optimization problems, and proved that it owns the $\mathcal{O}(1/k^2)$ convergence rate. It is worth noting that the convergence rate analysis of the  accelerated algorithms mentioned above was done  for the Lagrangian residual $\mathcal{L}(x^*, \lambda^*) -\mathcal{L}(x_k, \lambda_k)$. Recently, \cite{Xu2017} presented  an accelerated ALM  for solving the problem \eqref{ques}. By adapting parameters during the iterations, they proved that the objective residual and the feasibility violation both enjoy the $\mathcal{O}(1/k^2)$ convergence rate.   By applying Nesterov's technique, \cite{HeHF2021Inertial} proposed two  accelerated primal-dual algorithms, which  enjoy the $\mathcal{O}(1/k^2)$ convergence rate of the  objective residual and the feasibility violation. In terms of  scaling coefficients,  \cite{HeHF2021,Luo2021primal,Yan2020} proposed  accelerated ALM algorithms in different ways, and obtained the  fast convergence rates related to the scaling coefficients. By time discretization of a second-order dynamical system, \cite{Luo2021Accer} proposed new  accelerated primal-dual methods, and derived  the  $\mathcal{O}(1/k^2)$ convergence rate for the primal-dual gap, the feasibility violation and the objective residual  under the assumption that the objective function is  strongly convex. In the  case that $f$ has a Lipschitz continuous gradient, \cite{Bot_arx} proposed  fast ALM algorithms by time discretization of the dynamical system in \cite{Bot}. They proved  the  $\mathcal{O}(1/k^2)$ convergence rate of the primal-dual gap, the feasibility measure and  the objective residual, and also showed the convergence of the sequence of iterations. From the variational perspective, \cite{FazlyabACC} proposed  accelerated higher-order gradient methods by discretization of a second-order dual dynamical system and  exhibited the $\mathcal{O}(1/k^p)$ converges rate of the dual residual and the $\mathcal{O}(1/k^{p/2})$ rate of the feasibility violation under the assumption that the objective function is strong convex and  has a $(p-1)$--th Lipschitz gradient.

\subsection{Fast primal-dual  algorithm via dynamical system}
Dynamical system methods  have been recognized as efficient tools for solving  optimization problems in the literature. Dynamical systems can not only give more insights into the existing numerical methods for optimization problems but also  lead to other possible numerical algorithms by  time discretization, see e.g. \cite{ChenL2021,Liang2019,Jordan2018,Su2014,Wilson2021,Attouch2018MP,Luo2021Uni,Kia2015}. In this paper,  we will  propose a fast primal-dual  algorithm via  time discretization of the following primal-dual dynamical system
\begin{equation}\label{dy_pri_dual_a}
	\begin{cases}
		\ddot{x}(t)+\frac{\alpha}{t}\dot{x}(t)&= -\beta(t)(\nabla f(x(t))+A^T\lambda(t))+\epsilon(t),\\
		\dot{\lambda}(t)&= t\beta(t)(A(x(t)+ \frac{t}{\alpha-1}\dot{x}(t))-b),
	\end{cases}
\end{equation}
where $t\geq t_0>0$, $\alpha/t$ with $\alpha>1$ is a damping coefficient, $\beta:[t_0,+\infty)\to (0,+\infty)$ is a scaling coefficient, and  $\epsilon:[t_0,+\infty)\to \mathbb{R}^n$ is a perturbation coefficient. 

\cite{Su2014}  showed that the damping  ${\alpha}/{t}$ with $\alpha = 3$ in inertial dynamical systems  for unconstrained convex optimization problems can be understood  as the continuous limit of  Nesterov's accelerated technique  \cite{Nesterov1983}.  By considering the damping  ${\alpha}/{t}$ with $\alpha\geq 3$ and the scaling $\beta(t)$,  \cite{Attouch2019FP,Wibisono2016PANS} obtained the  $\mathcal{O}(1/t^2\beta(t))$ convergence rate of dynamical systems for solving  unconstrained convex optimization problems and also obtained the rate-matching inertial algorithms  by different time discretization schemes.   Recently, some researchers extended the dynamical systems in  \cite{Su2014,Attouch2019FP,Wibisono2016PANS}  to inertial primal-dual dynamical systems for solving the problem \eqref{ques}. See e.g. \cite{HeHF2020Siam,Attouch2021jt,Zeng2019Arxiv,Bot}. Very recently, by time discretization of the inertial dynamical system, \cite{Bot_arx} proposed  new ALM algorithms with  $\mathcal{O}(1/k^2)$ convergence rate. It is worth mentioning that for the dynamical systems mentioned above, inertial constructions are needed for both the primal variable  and  the dual variable. As a comparison, the inertial term is considered in the dynamical system \eqref{dy_pri_dual_a} only for the primal variable.

In the next, by time discretization of \eqref{dy_pri_dual_a},  we will propose an accelerated primal-dual algorithm. Rewrite \eqref{dy_pri_dual_a}  as 
\begin{equation}\label{dy_rewrit}
	\begin{cases}
	y(t) &= x(t) + \frac{t}{\alpha-1}\dot{x}(t),\\
	\frac{\alpha-1}{t}\dot{y}(t) &= -\beta(t)(\nabla f(x(t))+A^T\lambda(t))+\epsilon(t),\\
		\dot{\lambda}(t) &= t\beta(t)(Ay(t)-b).
	\end{cases}
\end{equation}
Set $t_k=k, x_k=x(t_k),\ y_k=y(t_k),\ \lambda_k = \lambda(t_k),\ \beta_k = \beta(t_k)$, $\epsilon_k=\epsilon(t_k)$. Take the following discretization scheme of  \eqref{dy_rewrit}  with a  nondifferentiable function $f$
\begin{subequations}\label{discrete}
\begin{numcases}{}
y_k = x_k+\frac{k-\theta}{\alpha-1}(x_{k}-x_{k-1}),\label{dis1}\\
\frac{\alpha-1}{k}(y_{k+1}-y_k) \nonumber\\
\qquad\quad \in -\beta_k(\partial f(x_{k+1})+A^T\lambda_{k+1})+\epsilon_k,\label{dis2}\\
\lambda_{k+1}-\lambda_k =k\beta_k(Ay_{k+1}-b),\label{dis3}
 \end{numcases}
\end{subequations}
where $\theta\in\mathbb{R}$.
Substitute  \eqref{dis1} and \eqref{dis3} into  \eqref{dis2}  to obtain
\begin{eqnarray}\label{optim}
&&\frac{k+\alpha-\theta}{k}\left(x_{k+1}-\left(x_k+\frac{k-\theta}{k+\alpha-\theta}(x_k-x_{k-1})\right)\right)\nonumber\\
&& \  \in -\beta_k(\partial f(x_{k+1})+A^T\lambda_{k})+\epsilon_k\\
&&\quad -k\beta_k^2A^T\left(A\left(x_{k+1}+\frac{k+1-\theta}{\alpha-1}(x_{k+1}-x_{k})\right)-b\right).\nonumber
\end{eqnarray}
Then, we can write  \eqref{discrete} as Algorithm \ref{al:al1}  for solving problem \eqref{ques}.
\begin{algorithm}[h]
\setstretch{0.6} 
        \caption{Fast  primal-dual algorithm}\label{al:al1}
        {\bf Initialization:} Choose $x_0\in\mathbb{R}^n,\ \lambda_0\in\mathbb{R}^m$, $\alpha> 1$,  $\theta\in\mathbb{R}$. Set $x_1 =x_0,\ \lambda_1=\lambda_0$. \\
       {\bf For} {$k = 1, 2,\cdots$}{\bf do}
   
        {\quad \bf Step 1:} Compute $\bar{x}_k = x_k+\frac{k-\theta}{k+\alpha-\theta}(x_k-x_{k-1})$.
        
       {\quad \bf Step 2:} Choose  $\beta_k>0$. Set 
       \[\vartheta_k =\frac{k(k+\alpha-\theta)\beta_k}{\alpha-1}, \] 
       \[\eta_k = \frac{k+1-\theta}{k+\alpha-\theta}Ax_k+ \frac{\alpha-1}{k+\alpha-\theta}b.\]
  \qquad  Update the primal variable 
  \begin{eqnarray*}
\qquad  	x_{k+1} &=& \mathop{\arg\min}_{x\in\mathbb{R}^{n}}\ \{f(x)+\frac{k+\alpha-\theta}{2k\beta_k}\|x-\bar{x}_k\|^2\\
  	&&+\frac{\vartheta_k}{2}\|Ax-\eta_k\|^2+\langle A^T{\lambda}_k-\frac{\epsilon_k}{\beta_k}, x\rangle\}.
  \end{eqnarray*} 
   {\quad \bf Step 3:} Compute
   \[ y_{k+1} =x_{k+1}+\frac{k+1-\theta}{\alpha-1}(x_{k+1}-x_k).\]
    \qquad  Update the dual variable
   \[\quad\lambda_{k+1} =\lambda_k+k\beta_k(Ay_{k+1}-b).\]
   {\quad \bf If}  A stopping condition is satisfied  {\bf then}

	 {\quad\qquad \bf Return} $(x_{k+1},\lambda_{k+1})$
	 
	{\quad \bf end}\\
{\bf end}
\end{algorithm}
The perturbation $\epsilon(t)$ in \eqref{dy_pri_dual_a} can be interpreted as a kind of disturbance, and here we adopt the terminology ``perturbation'' used by \cite{Attouch2018MP,HeHF2020Siam}. In Step 2 of Algorithm \ref{al:al1}, the perturbation $\epsilon_k$ means 
 \begin{eqnarray*}
\qquad  	x_{k+1} &\approx& \mathop{\arg\min}_{x\in\mathbb{R}^{n}}\ \{f(x)+\frac{k+\alpha-\theta}{2k\beta_k}\|x-\bar{x}_k\|^2\\
  	&&+\frac{\vartheta_k}{2}\|Ax-\eta_k\|^2+\langle A^T{\lambda}_k, x\rangle\}.
  \end{eqnarray*} 
The $x$-subproblem in Step 2 of Algorithm \ref{al:al1} has a special splitting structure and it can be efficiently solved by some classical splitting  methods  such as the proximal gradient method and its accelerated version FISTA (see e.g. \cite{Beck2009,Lin2020}).

In this paper, by constructing a  discrete  energy sequence and a continuous energy function,  we show fast convergence properties of Algorithm \ref{al:al1} and the dynamical system \eqref{dy_pri_dual_a}.  Our main contributions are summarized as follows:

 {\bf (a):  The discrete level:} By  time discretization of the dynamical system \eqref{dy_pri_dual_a}  with damping  $\alpha/t$ and  scaling $\beta(t)$, we obtain  Algorithm \ref{al:al1}. Under a suitable scaling  condition, we obtain the $\mathcal{O}(1/k^2\beta_k)$ convergence rate of  the objective residual and the feasibility violation, which can reach $\mathcal{O}(1/k^{\alpha-1})$ decay rate at the most. We  extend  \cite[Theorem 3.1 and Theorem 7.1]{Attouch2019FP} from the unconstrained optimization problem to  the linearly  constrained convex optimization problem \eqref{ques}. Compared with the accelerated  gradient  methods in \cite{Bot_arx} where $f$ is  convex and has a Lipschitz gradient, and \cite{FazlyabACC} where $f$ is strongly
convex, Algorithm \ref{al:al1}  requires  neither strong convexity  nor Lipschitz gradient assumption  on $f$. In the case $\alpha>3$, we show that Algorithm \ref{al:al1} can  achieve a rate faster than  $\mathcal{O}(1/k^{2})$ under a suitable scaling condition.
 
{\bf (b)  The continuous level:} For a better understanding of  the  acceleration scheme of Algorithm \ref{al:al1},  we consider the  primal-dual  dynamical system \eqref{dy_pri_dual_a} and show that it enjoys   convergence properties matching to that of  Algorithm \ref{al:al1}. To the best of our knowledge, the dynamical system \eqref{dy_pri_dual_a} is the first  Nesterov's inertial one involving  inertial term only for the primal variable for the linearly constrained optimization problem.  Our dynamical system \eqref{dy_pri_dual_a}  extends  the dynamical system in \cite{HeHF2021}, which  is linked to Polyak' heavy ball scheme, from the constant viscous damping to the vanishing damping  $ {\alpha}/{t}$.  Compared with a recent work by \cite{Attouch2021jt},  where  a general ADMM dynamical system involving inertial terms both for the primal and dual variables was considered for a separable linearly constrained optimization problem,  by a new result (Lemma \ref{le_A4}), we will prove that the convergence results of the dynamical system \eqref{dy_pri_dual_a} are better than the one in \cite[Theorem 1]{Attouch2021jt}. By  Lemma \ref{le_A4}, we also can improve the convergence of the objective residual and  the feasibility violation in  \cite[Theorem 1]{Attouch2021jt} from $\mathcal{O}(1/t^{{1}/{2\alpha_0}})$ to $\mathcal{O}(1/t^{{1}/{\alpha_0}})$. In the case  $\beta(t)\equiv\beta>0$,  the dynamical system \eqref{dy_pri_dual_a} enjoys the same convergence rate as the dynamical systems in \cite{Bot,Zeng2019Arxiv}, which involve inertial terms both for the primal and dual variables.

\subsection{Organization}
The paper is organized as follows: In Section 2,  we show the fast convergence properties of Algorithm \ref{al:al1} under a suitable scaling condition.    Section 3 is devoted to the study of convergence properties of the inertial primal-dual  dynamical system \eqref{dy_pri_dual_a}. The numerical experiments are given in  Section 4. Finally, we end the paper with a conclusion.

\section{Fast convergence analysis of Algorithm \ref{al:al1}}

Before presenting the convergence analysis, we first show that Algorithm \ref{al:al1} is equivalent to the time discretization scheme \eqref{discrete}.
\begin{proposition}\label{eq_le_20}
	Algorithm \ref{al:al1} is equivalent to the time discretization scheme \eqref{discrete}.
\end{proposition}
 \begin{pf}
 By using the optimality criterion, from Step 2 of  Algorithm \ref{al:al1}, we get
 \begin{eqnarray*}\label{eq_le21}
 	  0&\in& \partial f(x_{k+1})+\frac{k+\alpha-\theta}{k\beta_k}(x_{k+1}-\bar{x}_k)\\
 	  &&+A^T(\vartheta_k (Ax_{k+1}-\eta_k)+{\lambda}_k)-\frac{\epsilon_k}{\beta_k} ,
 \end{eqnarray*}
 which can be rewritten  as
  \begin{eqnarray}\label{eq_le22}
 	&&\frac{k+\alpha-\theta}{k}(x_{k+1}-\bar{x}_k) \in -\beta_k(\partial f(x_{k+1})\nonumber\\
 	&&\qquad+A^T(\vartheta_k (Ax_{k+1}-\eta_k)+{\lambda}_k))+\epsilon_k.
 \end{eqnarray}
It follows from Step 2 and Step 3 of  Algorithm \ref{al:al1} that
\begin{eqnarray}\label{eq_le23}
 	&&\vartheta_k(Ax_{k+1}-\eta_k) = \frac{k(k+\alpha-\theta)\beta_k}{\alpha-1}A x_{k+1}\nonumber\\
 	&&\qquad -\frac{k(k+1-\theta)\beta_k}{\alpha-1} Ax_{k}-k\beta_k b \\
 	&&\quad =  k\beta_k(A(x_{k+1}+\frac{k+1-\theta}{\alpha-1}(x_{k+1}-x_{k}))-b).\nonumber
 \end{eqnarray}
As a consequence of \eqref{eq_le22}, \eqref{eq_le23} and Step 1,  the equation   \eqref{optim} holds. By comparing Algorithm \ref{al:al1} and \eqref{discrete}, the sequence $\{(x_k,y_k,\lambda_k)\}_{k\geq 1}$  generated by Algorithm \ref{al:al1}  satisfies  \eqref{discrete}. Since the calculation process from above is reversible, from \eqref{discrete}, we also can obtain Algorithm \ref{al:al1}.
 \end{pf}

\subsection{Convergence analysis for fast primal-dual algorithm}
Before discussing  the  convergence properties of  Algorithm \ref{al:al1}, we first recall the   equality
\begin{equation}\label{eq_know1}
	 \frac{1}{2}\|x\|^2-\frac{1}{2}\|y\|^2=\langle x,x-y\rangle-\frac{1}{2}\|x-y\|^2
\end{equation}
for any $x,y,z \in\mathbb{R}^n$, which will be used repeatedly.
\begin{lemma}\label{le_ex08}
	Let  $\{(x_k,y_k,\lambda_k)\}_{k\geq 1}$ be the sequence generated by Algorithm \ref{al:al1} and $(x^*,\lambda^*)\in\Omega$. 
Define the  energy sequence
\begin{equation}\label{eq_th_aa}
	\mathcal{E}^{\epsilon}_k= \mathcal{E}_k-\sum^k_{j=1} \langle (\alpha-1)(y_{j}-x^*), (j-1)\epsilon_{j-1}\rangle
\end{equation}
with
\begin{eqnarray}\label{eq_th2_1}
	\mathcal{E}_k& =&k(k+1-\theta)\beta_k(\mathcal{L}(x_k,\lambda^*)-\mathcal{L}(x^*,\lambda^*))\nonumber \\
	&&+ \frac{1}{2}\|(\alpha-1)(y_k-x^*)\|^2+ \frac{\alpha-1}{2}\|\lambda_k-\lambda^*\|^2.
\end{eqnarray}
Then, for any $k\geq \max\{1,\theta-1\}$:
\begin{eqnarray}\label{eq_th2_12}
	&& \mathcal{E}^{\epsilon}_{k+1}-\mathcal{E}^{\epsilon}_k\nonumber\\
	&&\qquad \leq ((k+1)(k+2-\theta)\beta_{k+1}-k\left(k+\alpha-\theta)\beta_k\right)\nonumber\\
	&&\quad\qquad \cdot(\mathcal{L}(x_{k+1},\lambda^*)-\mathcal{L}(x^*,\lambda^*)).
\end{eqnarray}
\end{lemma}
\begin{pf}
By the definition of $\mathcal{L}$, we have
\begin{equation*}\label{eq_th2_3}
	 \partial_x \mathcal{L}(x,\lambda) = \partial f(x)+ A^T\lambda.
\end{equation*}
This together with   \eqref{dis2} implies
\begin{eqnarray*}\label{eq_th2_4}
	&& (\alpha-1)(y_{k+1}-y_k)\nonumber \\
	&&\in k(-\beta_k(\partial f(x_{k+1})+ A^T\lambda_{k+1})+\epsilon_k)\\
	&& = k(-\beta_k(\partial f(x_{k+1})+ A^T\lambda^*)-\beta_k A^T(\lambda_{k+1}-\lambda^*)+\epsilon_k) \nonumber\\
	&&= -k\beta_k \partial_x \mathcal{L}(x_{k+1},\lambda^*)- k\beta_kA^T(\lambda_{k+1}-\lambda^*)+k\epsilon_k.\nonumber
\end{eqnarray*}
Denote
\begin{eqnarray}\label{eq_th2_5}
	\quad \xi_k &:=& -\frac{\alpha-1}{k\beta_k}(y_{k+1}-y_k)-A^T(\lambda_{k+1}-\lambda^*)+\frac{\epsilon_k}{\beta_k} \nonumber\\
	&\in &\partial_x \mathcal{L}(x_{k+1},\lambda^*).
\end{eqnarray}
Combining \eqref{eq_know1},  \eqref{eq_th2_5} and $\alpha>1$, we have
\begin{eqnarray}\label{eq_th2_6}
	&&\frac{1}{2}\|(\alpha-1)(y_{k+1}-x^*)\|^2- \frac{1}{2}\|(\alpha-1)(y_k-x^*)\|^2\nonumber \\
	&& = \langle (\alpha-1)(y_{k+1}-x^*), (\alpha-1)(y_{k+1}- y_{k})\rangle\nonumber\\
	&&\quad -\frac{(\alpha-1)^2}{2}\| y_{k+1}- y_{k}\|^2\\
	 &&\leq -k\beta_k \langle (\alpha-1)(y_{k+1}-x^*)  ,\xi_k \rangle
	\nonumber \\
	 &&\quad -  k\beta_k\langle(\alpha-1)(y_{k+1}-x^*), A^T(\lambda_{k+1}-\lambda^*)\rangle \nonumber\\
	 &&\quad +\langle (\alpha-1)(y_{k+1}-x^*),k\epsilon_k\rangle.\nonumber
\end{eqnarray}
From assumption we have $k+1-\theta\geq 0$ and $\alpha> 1$. By  \eqref{eq_th2_5} and Step 3, we get
\begin{eqnarray}\label{eq_th2_7}
	&&\langle (\alpha-1)(y_{k+1}-x^*), \xi_k\rangle\nonumber \\
	&&\  =(\alpha-1)\langle x_{k+1}-x^*,  \xi_k\rangle+(k+1-\theta)\langle x_{k+1}-x_k, \xi_k \rangle\nonumber \\
	&&\  \geq(\alpha-1) (\mathcal{L}(x_{k+1},\lambda^*)-\mathcal{L}(x^{*},\lambda^*))\\
	&&\quad+(k+1-\theta)(\mathcal{L}(x_{k+1},\lambda^*)-\mathcal{L}(x_{k},\lambda^*)),\nonumber
\end{eqnarray}
where the inequality follows from the convexity of $\mathcal{L}(\cdot,\lambda^*)$. 
By  Step 3,  $Ax^*= b$, and \eqref{eq_know1}, we get
\begin{eqnarray}\label{eq_th2_10}
	&&\frac{1}{2}\|\lambda_{k+1}-\lambda^*\|^2-\frac{1}{2}\|\lambda_{k}-\lambda^*\|^2\nonumber\\
	&&\quad= \langle \lambda_{k+1}-\lambda^*, \lambda_{k+1}-\lambda_k\rangle -\frac{1}{2}\|\lambda_{k+1}-\lambda_k\|^2 \\
	&&\quad\leq  \langle \lambda_{k+1}-\lambda^*,  k\beta_k A(y_{k+1}-x^*)\rangle.\nonumber
\end{eqnarray}
It follows from \eqref{eq_th_aa} and \eqref{eq_th2_1}  that 
\begin{eqnarray*}
	&&\mathcal{E}^{\epsilon}_{k+1}-\mathcal{E}^{\epsilon}_k = \mathcal{E}_{k+1}-\mathcal{E}_k-\langle(\alpha-1)(y_{k+1}-x^*),k\epsilon_k\rangle\nonumber\\
	&& \  \leq ((k+1)(k+2-\theta)\beta_{k+1}-k(k+\alpha-\theta)\beta_k)\\
	&&\quad \cdot(\mathcal{L}(x_{k+1},\lambda^*)-\mathcal{L}(x^*,\lambda^*))\nonumber\\
	&&\quad-  k\beta_k\langle(\alpha-1)(y_{k+1}-x^*), A^T(\lambda_{k+1}-\lambda^*)\rangle \\
	&&\quad+\frac{\alpha-1}{2}(\|\lambda_{k+1}-\lambda^*\|^2-\|\lambda_{k}-\lambda^*\|^2)\nonumber\\
	&& \ \leq ((k+1)(k+2-\theta)\beta_{k+1}-k(k+\alpha-\theta)\beta_k)\\
	&&\quad \cdot(\mathcal{L}(x_{k+1},\lambda^*)-\mathcal{L}(x^*,\lambda^*)),\nonumber
\end{eqnarray*}
where the first inequality follows from \eqref{eq_th2_6} and \eqref{eq_th2_7}, and the last inequality follows from \eqref{eq_th2_10}. This yields the desired result.
\end{pf}

To derive the fast convergence rates, we need the following scaling condition: there exist $k_1\geq\max\{2,\theta\}$ such that
\begin{equation}\label{eq_assum_beta}
	\beta_{k+1}\leq \frac{k(k+\alpha-\theta)}{(k+1)(k+2-\theta)}\beta_k, \quad\forall k\geq k_1-1.
\end{equation}
 Now, we start to discuss the fast convergence properties of Algorithm \ref{al:al1} by the Lyapunov analysis approach.
\begin{theorem}\label{th_al1}
Let  $\{(x_k,y_k,\lambda_k)\}_{k\geq 1}$ be the sequence generated by Algorithm \ref{al:al1} and $(x^*,\lambda^*)\in\Omega$. Assume that the condition \eqref{eq_assum_beta} holds and 
\[\sum^{+\infty}_{k=1} k\|\epsilon_k\| <+\infty,\quad \lim_{k\to+\infty}k^2\beta_k=+\infty.\]
Then, the sequence $\{(y_k,\lambda_k)\}_{k\geq k_1}$ is bounded,
\[\|Ax_k-b\|=\mathcal{O}\left(\frac{1}{k^2\beta_k}\right),\] 
and
	\[ 	|f(x_k)-f(x^*)|=\mathcal{O}\left(\frac{1}{k^2\beta_k}\right).\]
\end{theorem}

\begin{pf}
From assumptions, we can get $\mathcal{L}(x_{k+1},\lambda^*)-\mathcal{L}(x^*,\lambda^*)\geq 0$ and $(k+1)(k+2-\theta)\beta_{k+1}-k(k+\alpha-\theta)\beta_k\leq 0$ for any $k\geq k_1$.     Then, for any $k\geq k_1$, $\mathcal{E}_{k}\geq  0$, and  from Lemma \ref{le_ex08} we have
\begin{equation*}\label{eq_th3_12}
	\mathcal{E}^{\epsilon}_{k}\leq \mathcal{E}^{\epsilon}_{k_1}, \qquad \forall k\geq k_1.
\end{equation*}
By \eqref{eq_th_aa} and \eqref{eq_th2_1}, we have
\begin{eqnarray}\label{eq_th3_13}
	&& \frac{1}{2}\|(\alpha-1)(y_k-x^*)\|^2\leq\mathcal{E}_{k}\nonumber\\
	&&\qquad =\mathcal{E}^{\epsilon}_{k}+\sum^k_{j=1}\langle (\alpha-1)(y_j-x^*),(j-1)\epsilon_{j-1}\rangle\nonumber \\
	&&\qquad \leq \mathcal{E}^{\epsilon}_{k_1}+\sum^k_{j=1}\langle (\alpha-1)(y_j-x^*),(j-1)\epsilon_{j-1}\rangle \\
	&&\qquad= \mathcal{E}_{k_1}+\sum^k_{j=k_1+1}\langle (\alpha-1)(y_j-x^*),(j-1)\epsilon_{j-1}\rangle\nonumber \\
	&&\qquad \leq \mathcal{E}_{k_1}+\sum^{k}_{j=k_1}(j-1)\|(\alpha-1)(y_j-x^*)\|\cdot\|\epsilon_{j-1}\|\nonumber
\end{eqnarray}
for any $k\geq k_1$. 
Note that $\sum^{+\infty}_{k=1} k\|\epsilon_k\| <+\infty$. Applying Lemma \ref{le_disc_per} with $a_k=\|(\alpha-1)(y_{k+k_1-1}-x^*)\|$, we get 
\begin{equation*}\label{eq_th3_14}
	\sup_{k\geq k_1}\|(\alpha-1)(y_k-x^*)\|=	\sup_{k\geq 1}a_k<+\infty.
\end{equation*}
This  together with \eqref{eq_th3_13} yields
\begin{equation*}\label{eq_th3_15}
	\sup_{k\geq k_1}\mathcal{E}_{k}\leq \mathcal{E}_{k_1}+ \sup_{k\geq k_1}\|(\alpha-1)(y_k-x^*)\|\cdot\sum^{+\infty}_{j=1}j\|\epsilon_j\|<+\infty.
\end{equation*}
Thus, the energy sequence $\{\mathcal{E}_{k}\}_{k\geq k_1}$ is  bounded.  By \eqref{eq_th2_1}, $\{\|y_k-x^*\|\}_{k \geq k_1}$ and  $\{\|\lambda_k-\lambda^*\|\}_{k\geq k_1}$  are  bounded and 
\begin{equation}\label{eq_th3_16}
	\sup_{k\geq k_1}k(k+1-\theta)\beta_k(\mathcal{L}(x_k,\lambda^*)-\mathcal{L}(x^*,\lambda^*))<+\infty. 
\end{equation}
As a result,  the sequence $\{(y_k,\lambda_k)\}_{k\geq k_1}$ is bounded. From \eqref{eq_th3_16} we obtain 
\begin{equation}\label{eq_L_kbeta}
	\mathcal{L}(x_k,\lambda^*)-\mathcal{L}(x^*,\lambda^*)=\mathcal{O}\left({1}/{k^2\beta_k}\right).
\end{equation}

For notation simplicity, denote 
\begin{equation}\label{eq_hx1}
g_k :=\frac{(k+\alpha-\theta-1)(k-1)\beta_{k-1}}{\alpha-1}(Ax_{k}-b).
\end{equation}
It follows from  Step 3 that
\begin{eqnarray}\label{eq_g_comp}
	&&\lambda_{k+1}-\lambda_{k_1}=\sum^{k}_{j=k_1}(\lambda_{j+1}-\lambda_j)\nonumber \\
	&&=\sum^{k}_{j=k_1}j\beta_j(Ay_{j+1}-b)\\
	&&=  \sum^{k}_{j=k_1} j\beta_j\left((Ax_{j+1}-b)+\frac{j+1-\theta}{\alpha-1}A(x_{j+1}-x_j)\right)\nonumber \\
	&& =\sum^{k}_{j=k_1}\left(g_{j+1}-\frac{(j+1-\theta)j\beta_j}{(j+\alpha-\theta-1)(j-1)\beta_{j-1}}g_j\right)\nonumber
\end{eqnarray}
for any $k\geq k_1$.
Let 
\[a_k = 1-\frac{(k+1-\theta)k\beta_k}{(k+\alpha-\theta-1)(k-1)\beta_{k-1}}, \quad \forall k\geq k_1.\]
From \eqref{eq_g_comp},  we have
\begin{eqnarray*}
	\lambda_{k+1}-\lambda_{k_1}&=&\sum^{k}_{j=k_1}(g_{j+1}- g_j)+ \sum^{k}_{j=k_1}a_jg_j\nonumber\\
	&=& g_{k+1}-g_{k_1} + \sum^{k}_{j=k_1}a_jg_j.
\end{eqnarray*}
This together with the boundedness of $\{\lambda_k\}_{k\geq k_1}$ yields
\[\left\|g_{k+1} +\sum^{k}_{j=k_1}a_jg_j\right\|\leq C, \quad \forall k\geq k_1,\]
where 
\[C = \sup_{k\geq k_1} \|\lambda_{k+1}-\lambda_{k_1}\| +\|g_{k_1}\|<+\infty.\]
From \eqref{eq_assum_beta}, we can get 
\[0\leq a_k <1,\quad \forall k\geq k_1.\]
Applying  Lemma \ref{le_A2}, we have 
\[\sup_{k\geq k_1}\|g_k\|<+\infty,\]
which together with \eqref{eq_assum_beta} and \eqref{eq_hx1} implies
\[\|Ax_k-b\| = \mathcal{O}\left(\frac{1}{k^2\beta_k}\right).\]
It follows from \eqref{eq_L_kbeta} that
\begin{eqnarray*}
&&	|f(x_k)-f(x^*)|\\
&&\qquad\leq  \mathcal{L}(x_k,\lambda^*)-\mathcal{L}(x^*,\lambda^*)+\|\lambda^*\|\|Ax_k-b\| \\
&&	\qquad= \mathcal{O}\left(\frac{1}{k^2\beta_k}\right).
\end{eqnarray*}
\end{pf}

\begin{remark}
The assumption \eqref{eq_assum_beta} is just the assumption $(H_{\beta,\theta})$ appears in \cite{Attouch2019FP} for  convergence rate analysis of inertial proximal algorithms for the unconstrained optimization problems, and  Theorem \ref{th_al1} can be viewed as an extension of \cite[Theorem 3.1 and Theorem 7.1]{Attouch2019FP} from the  unconstrained case to the problem \eqref{ques}.
\end{remark}

From Theorem \ref{th_al1}, we can obtain the best decay rate when the condition \eqref{eq_assum_beta} holds with equality, such that 
\begin{equation}\label{eq_assum_beta_eq}
	\beta_{k+1}= \frac{k(k+\alpha-\theta)}{(k+1)(k+2-\theta)}\beta_k, \quad\forall k\geq k_1-1
\end{equation}
with $k_1\geq \max\{2,\theta\}$.
\begin{corollary}\label{corr1}
Suppose the assumptions of Theorem \ref{th_al1} hold and that $\{\beta_k\}_{k\geq 1}$ satisfies \eqref{eq_assum_beta_eq}.
Let  $\{(x_k,y_k,\lambda_k)\}_{k\geq 1}$ be the sequence generated by Algorithm \ref{al:al1} and $(x^*,\lambda^*)\in\Omega$. Then, 
\[\|Ax_k-b\|=\mathcal{O}\left(\frac{1}{k^{\alpha-1}}\right),\] 
	\[ |f(x_k)-f(x^*)|=\mathcal{O}\left(\frac{1}{k^{\alpha-1}}\right).\]
\end{corollary}

\begin{pf}
	From \eqref{eq_assum_beta_eq},  we have
	\[(k+1)(k+2-\theta)\beta_{k+1} = (1+\frac{\alpha-1}{k+1-\theta})k(k+1-\theta)\beta_k\]
for all $k\geq k_1$. 
Let 
\begin{equation}\label{eq_new_315}
	\gamma_k = (k+k_1-1)(k+k_1-\theta)\beta_{k+k_1-1},\quad \forall k\geq 1.
\end{equation}
Then, 
\[\gamma_{k+1} =(1+\frac{\alpha-1}{k+k_1-\theta})\gamma_{k},\quad \forall k\geq 1.\]
By Lemma \ref{le_new}, there exists $\mu_1>0$ and $\mu_2>0$ such that 
\[ \mu_1 k^{\alpha-1}\leq  \gamma_k\leq \mu_2 k^{\alpha-1}.\]  This together with \eqref{eq_new_315} and  Theorem \ref{th_al1} yields the desired result. 
\end{pf}

\begin{remark}
Under the assumption that $f$ is strongly convex and has an $(\alpha-2)$--th Lipschitz gradient,    \cite{FazlyabACC} proposed an $\mathcal{O}(1/k^{\alpha-1})$ convergence rate algorithm for  the problem \eqref{ques}. As a comparison,   Algorithm \ref{al:al1} can enjoy the $\mathcal{O}(1/k^{\alpha-1})$ decay rate under the merely convexity assumption. 
\end{remark}

 Let $Id$ be the identity matrix and $\mathbb{S}_+(n)$ be the set of all positive semidefinite matrixes in $\mathbb{R}^{n\times n}$. Denote $\|x\|^2_M =x^TMx$ for any $x\in\mathbb{R}^{n}$ and $M\in\mathbb{S}_+(n)$. For any $M_1, M_2 \in \mathbb{S}_+(n)$, denote
 \[M_1 \succcurlyeq M_2 \Longleftrightarrow \|x\|_{M_1}\geq \|x\|_{M_2},\qquad\forall x\in\mathbb{R	}^n .\]
  It is easy to verify that  for any $x,y \in\mathbb{R}^{n}$,  $M\in\mathbb{S}_+(n)$,
 \begin{equation*}\label{eq_know3}
 	 \frac{1}{2}\|x\|_M^2-\frac{1}{2}\|y\|_M^2 = \langle x, M(x-y)\rangle-\frac{1}{2}\|x-y\|_M^2. 
 \end{equation*}
Then,   we can replace the subproblem of step 2 with
\begin{eqnarray}\label{subproblem_M}
	x_{k+1} &= & \mathop{\arg\min}_{x\in\mathbb{R}^{n}}\ \{ f(x)+\frac{k+\alpha-\theta}{2k\beta_k}\|x-\bar{x}_k\|_M^2\nonumber\\
	&& +\frac{\vartheta_k}{2}\|Ax-\eta_k\|^2+\langle A^T{\lambda}_k-\frac{\epsilon_k}{\beta_k}, x\rangle\},
\end{eqnarray}
where $M\succcurlyeq \kappa Id$ for some $\kappa > 0$. Redefine \eqref{eq_th2_1} as
\begin{eqnarray*}
	\mathcal{E}_k&=&k(k+1-\theta)\beta_k(\mathcal{L}(x_k,\lambda^*)-\mathcal{L}(x^*,\lambda^*)	)\nonumber \\
	&&+ \frac{1}{2}\|(\alpha-1)(y_k-x^*)\|_M^2+ \frac{\alpha-1}{2}\|\lambda_k-\lambda^*\|^2.
\end{eqnarray*}
 Through the arguments similar to the ones in Theorem \ref{th_al1}, we can get the same convergence  results. In particular, when the perturbation $\epsilon_k\equiv 0$, which means that the subproblems are solved with exact or high precision, we can take $\kappa=0$.

\section{Convergence  properties of   inertial primal-dual  dynamical  system}

In this section, for a better understanding of the acceleration scheme of  Algorithm \ref{al:al1}, we will investigate the convergence properties of the dynamical system \eqref{dy_pri_dual_a}. When $\nabla f$ is globally Lipschitz continuous, $\beta(t)$ is a continuous differentiable function, through the  Cauchy–Lipschitz theorem \cite[Proposition 6.2.1]{Haraux} and  the similar discussions in \cite[Theorem 5]{Attouch2021jt}, we can prove that \eqref{dy_pri_dual_a} has a unique strong global $C^2$ solution.
 In what follows, we always assume that  $f$ is a proper,  convex and differentiable function and that  \eqref{dy_pri_dual_a} admits a  global solution.

\begin{theorem}\label{th_A1}
	 Assume that  $\alpha>1$,  $\beta:[t_0,+\infty)\to(0,+\infty)$ is a continuous differentiable function satisfying
\begin{equation}\label{ass_22}
	t\dot{\beta}(t) \leq (\alpha-3) \beta(t), \quad  \lim_{t\to+\infty}t^2\beta(t) = +\infty,
\end{equation}
and  $\epsilon:[t_0,+\infty)\to\mathbb{R}^n$ is an  integrable function satisfying
\begin{equation*}
\int^{+\infty}_{t_0}t\|\epsilon(t)\|dt<+\infty.
\end{equation*}
Let $(x(t), \lambda (t))$ be a global solution of the dynamical system \eqref{dy_pri_dual_a} and $(x^*,\lambda^*)\in\Omega$.  Then, 
\[\|Ax(t)-b\| =\mathcal{O}\left(\frac{1}{t^2\beta(t)}\right),\]
\[ |f(x(t))-f(x^*)| =\mathcal{O}\left(\frac{1}{t^2\beta(t)}\right).\]
\end{theorem}
\begin{pf}
 Define the energy  function  $\mathcal{E}^{\epsilon}:[t_0,+\infty)\to\mathbb{R}$ as
\begin{eqnarray*}\label{eq_th5_1}
	\mathcal{E}^{\epsilon}(t) =\mathcal{E}(t)- \int^t_{t_0}\langle (\alpha-1)(y(s)-x^*),s\epsilon(s)\rangle ds,
\end{eqnarray*}
where  $\mathcal{E}(t) = \mathcal{E}_0(t)+\mathcal{E}_1(t) $ with
\begin{equation}\label{ener_conti}
	\begin{cases}
		\mathcal{E}_0(t)&=t^2\beta(t)(\mathcal{L}(x(t),\lambda^*)-\mathcal{L}(x^*,\lambda^*)	),\\
		\mathcal{E}_1(t) &= \frac{1}{2}\|(\alpha-1)(y(t)-x^*)\|^2+\frac{\alpha-1}{2}\|\lambda(t)-\lambda^*\|^2,
	\end{cases}
\end{equation}
and $y(t)$ is defined by \eqref{dy_rewrit}.
By the classical differential calculations, \eqref{saddle_point} and \eqref{dy_rewrit}, we have
\begin{eqnarray*}
	\dot{\mathcal{E}}_0(t)&=&t^2\beta(t)\langle\nabla f(x(t))+A^T\lambda^*,\dot{x}(t)\rangle\\
	&& +(2t\beta(t)+t^2\dot{\beta}(t))(\mathcal{L}(x(t),\lambda^*)-\mathcal{L}(x^*,\lambda^*))
\end{eqnarray*}
and
\begin{eqnarray*}
	\dot{\mathcal{E}}_1(t)&=&\langle(\alpha-1)(y(t)-x^*),(\alpha-1)\dot{y}(t)\rangle\\
	&& +(\alpha-1)\langle \lambda(t)-\lambda^*,\dot{\lambda}(t)\rangle\\
	&=& -(\alpha-1)t\beta(t)\langle y(t)-x^*,\nabla f(x(t))+A^T\lambda(t)\rangle\\
	&&+\langle (\alpha-1)(y(t)-x^*),t\epsilon(t)\rangle\\
	 &&+(\alpha-1)t\beta(t)\langle \lambda(t)-\lambda^*,A(y(t)-x^*)\rangle \\
	&=&-(\alpha-1)t\beta(t)\langle x(t)-x^*, \nabla f(x(t))+A^T\lambda^*\rangle\\
	&& - t^2\beta(t)\langle \dot{x}(t), \nabla f(x(t))+A^T\lambda^*\rangle\\
	&&+\langle (\alpha-1)(y(t)-x^*),t\epsilon(t)\rangle.
\end{eqnarray*}
By computation, we  get
\begin{eqnarray}\label{eq_th5_2}
		&&\dot{\mathcal{E}^{\epsilon}}(t)= \dot{\mathcal{E}}_0(t)+\dot{\mathcal{E}}_1(t)-\langle (\alpha-1)(y(t)-x^*),t\epsilon(t)\rangle  \nonumber\\
	&&= (\alpha-1)t\beta(t)(f(x(t))-f(x^*)-\langle x(t)-x^*, \nabla f(x(t))\rangle)\nonumber \\
	&&\quad+ t(t\dot{\beta}(t)-  (\alpha-3)\beta(t))(\mathcal{L}(x(t),\lambda^*)-\mathcal{L}(x^*,\lambda^*)) \\
	&&\leq t(t\dot{\beta}(t)-  (\alpha-3)\beta(t))(\mathcal{L}(x(t),\lambda^*)-\mathcal{L}(x^*,\lambda^*)), \nonumber
\end{eqnarray}
where the inequality follows from the convexity of $f$. Since $(x^*,\lambda^*)\in\Omega$,  it is easy to verify that $\mathcal{L}(x(t),\lambda^*)-\mathcal{L}(x^*,\lambda^*)\geq 0$ and  $\mathcal{E}(t)\geq 0$. By assumptions and \eqref{eq_th5_2}, we get that $\dot{\mathcal{E}}^{\epsilon}(t)\leq 0$. As a result,
\begin{equation*}\label{eq_th5_3}
	\mathcal{E}^{\epsilon}(t)\leq  \mathcal{E}^{\epsilon}(t_0)=\mathcal{E}(t_0),\quad \forall t\in[ t_0,+\infty).
\end{equation*}
By  the definitions of $\mathcal{E}(t)$  and $\mathcal{E}^{\epsilon}(t)$,  and using the Cauchy-Schwarz inequality, for any $t\in[ t_0,+\infty)$ we have
\begin{eqnarray}\label{eq_hx_118}
	&&\frac{1}{2}\|(\alpha-1)(y(t)-x^*)\|^2\leq \mathcal{E}(t)\nonumber\\
	&&\ =\mathcal{E}^{\epsilon}(t)+\int^t_{t_0}\langle (\alpha-1)(y(s)-x^*),s\epsilon(s)\rangle ds,\\
	&&\  \leq \mathcal{E}(t_0)+\int^t_{t_0}\|(\alpha-1)(y(s)-x^*)\|\cdot s\|\epsilon(s)\| ds.\nonumber
\end{eqnarray}
Apply Lemma \ref{le_dy_per} with $\mu(t) =\|(\alpha-1)(y(t)-x^*)\| $ to get
\begin{eqnarray*}
&&\sup_{t\in[t_0,+\infty)}\|(\alpha-1)(y(t)-x^*)\|\\
 &&\quad \leq  \sqrt{2\mathcal{E}(t_0)}+\int^{+\infty}_{t_0}s\|\epsilon(s)\| ds<+\infty.
\end{eqnarray*}
This together with \eqref{eq_hx_118} implies
\begin{eqnarray*}
	&&\sup_{t\in[t_0,+\infty)}\mathcal{E}(t) \leq \mathcal{E}(t_0)\\
	&&\qquad+\sup_{t\in[t_0,+\infty)}\|(\alpha-1)(y(t)-x^*)\|\cdot\int^{+\infty}_{t_0}s\|\epsilon(s)\| ds \\
	&&\quad < +\infty.
\end{eqnarray*}
So, ${\mathcal{E}}(t)$ is  bounded, and then from \eqref{ener_conti} we get
 \begin{equation}\label{eq_thA_7}
	\mathcal{L}(x(t),\lambda^*)-\mathcal{L}(x^*,\lambda^*) =\mathcal{O}\left(\frac{1}{t^2\beta(t)}\right),
\end{equation}
and  $\lambda(t)$ is bounded on $[t_0,+\infty)$.

 By the partial integration, we can compute
 \begin{eqnarray*}
 	&&\int^t_{t_0} s^{2}\beta(s)A\dot{x}(s)ds = \int^t_{t_0} s^{2}\beta(s)d(Ax(s)-b)\\
 	&& \quad = t^2\beta(t)(Ax(t)-b)-t_0^2\beta(t_0)(Ax(t_0)-b) \\
	&&\qquad-\int^t_{t_0} (2s\beta(s)+s^2\dot{\beta}(s))(Ax(s)-b)ds.
 \end{eqnarray*} 
Then, from the second equation of \eqref{dy_pri_dual_a} we have
	\begin{eqnarray}\label{he_102_1}
	&&(\alpha-1)(\lambda(t) - \lambda(t_0))=  \int^t_{t_0}(\alpha-1)\dot{\lambda}(s)ds\nonumber \\
	&&\ = (\alpha-1)\int^t_{t_0} s\beta(s)A(x(s)-b)ds+ \int^t_{t_0}s^2\beta(s)A\dot{x}(s)ds\nonumber \\
	&&\  = t^2\beta(t)(Ax(t)-b)-t_0^2\beta(t_0)(Ax(t_0)-b)\\
	&&\quad+\int^t_{t_0}a(s)s^2\beta(s)(Ax(s)-b)ds\nonumber,
	\end{eqnarray}
where
\[ a(s) = \frac{\alpha-3}{s}-\frac{\dot{\beta}(s)}{\beta(s)}.\] From \eqref{ass_22} and the boundedness of   $\lambda(t)$, we get $a(t)\geq 0$ and
\begin{equation}\label{he_102_2}
	\left\| t^2\beta(t)(Ax(t)-b)+\int^t_{t_0}a(s)s^2\beta(s)(Ax(s)-b)ds\right\|\leq C
\end{equation}
for all $t\geq t_0$,
where
\begin{eqnarray*}
	C &=&(\alpha-1)\sup_{t\geq t_0}\|\lambda(t) - \lambda(t_0)\|+\|t_0^{2}\beta(t_0)(Ax(t_0)-b)\|\\
	&<&+\infty. 
\end{eqnarray*}
Now, applying Lemma \ref{le_A4} with $g(t) =t^2\beta(t)(Ax(t)-b)$, we obtain 
\begin{equation*}
	\sup_{t\geq t_0}\|t^2\beta(t)(Ax(t)-b)\|=\sup_{t\geq t_0}\|g(t)\|<+\infty,
\end{equation*}
which is 
\[\|Ax(t)-b\| =\mathcal{O}\left(\frac{1}{t^{2}\beta(t)}\right).\] 
This together with \eqref{eq_thA_7} implies
\begin{eqnarray*}
	&&|f(x(t))-f(x^*)|\\
	&&\quad\leq  \mathcal{L}(x(t),\lambda^*)-\mathcal{L}(x^*,\lambda^*)+\|\lambda^*\|\|Ax(t)-b\| \\
	&&\quad= \mathcal{O}\left(\frac{1}{t^{2}\beta(t)}\right).
\end{eqnarray*}
\end{pf}

\begin{remark}
	Theorem  \ref{th_A1} generalizes  \cite[Theorem A.1]{Attouch2019FP} from  the unconstrained optimization problem  to the  problem \eqref{ques}. Taking $t\dot{\beta}(t) = (\alpha-3) \beta(t)$,  then $\beta(t)=\frac{\beta(t_0)}{t_0^{\alpha-3}} t^{\alpha-3}$. From Theorem \ref{th_A1} we obtain the best $\mathcal{O}(1/t^{\alpha-1})$ decay rate  of the objective residual  and  the feasibility violation.  By  contrast, under  the strong  convexity assumption of $f$, \cite{FazlyabACC} only obtained the $\mathcal{O}(1/t^{\alpha-1})$ convergence rate of the dual residual and the $\mathcal{O}(1/t^{\frac{\alpha-1}{2}})$ convergence rate of the  feasibility violation for their dual dynamical system with  ${\alpha}/{t}$ damping. 
\end{remark}

\begin{remark}
As a comparison with the results on dynamical systems in \cite{Bot,Attouch2021jt}, we use a different method (Lemma \ref{le_A4}) to prove the fast convergence results of  \eqref{dy_pri_dual_a}.  By our method, we can simplify the proof process of \cite[Theorem 3.4]{Bot}, and also can  improve the convergence rate results  of the objective residual and the  feasibility violation in  \cite[Theorem 1]{Attouch2021jt} from $\mathcal{O}(1/t^{{1}/{2\alpha_0}})$ to $\mathcal{O}(1/t^{{1}/{\alpha_0}})$.
\end{remark}

\begin{remark}
The growth condition \eqref{eq_assum_beta} can be  rewritten as  
\[\beta_{k+1}-\beta_{k} \leq \frac{(\alpha-3)k+\theta-2}{(k+1)(k+2-\theta)}\beta_k.\]
This can be viewed as a discretized version of $\dot{\beta}(t)\leq \frac{\alpha-3}{t}\beta(t)$, which is \eqref{ass_22}. From Corollary \ref{corr1}, we  know that the scaling $\beta_k$ with \eqref{eq_assum_beta_eq} has the same order as $k^{\alpha-3}$, so it is the same order as the continuous function $\mu t^{\alpha-3}$ with $\mu>0$. In this sense, Theorem  \ref{th_A1} provides a dynamical interpretation of the fast convergence properties of Algorithm \ref{al:al1}.
\end{remark}

\section{Numerical experiments}
In this section, we test  Algorithm \ref{al:al1} on solving  the linearly constrained $\ell_1-\ell_2$ minimization problem. The numerical results demonstrate the validity and superior performance of our algorithm over some existing accelerated algorithms.

Consider the $\ell_1-\ell_2$ minimization problem
\begin{equation*}
	\min_x\ \|x\|_1 +\frac{\delta}{2}\|x\|^2_2 \quad s.t.\ Ax =  b,
\end{equation*}
where $A\in \mathbb{R}^{m\times n}$ and $b\in\mathbb{R}^m$. Set $m=1500,\ n=3000$ and $\delta =0.1$. Generate $A$ by the  standard Gaussian distribution and the original solution (signal) ${x}^*\in\mathbb{R}^n$ by the  Gaussian distribution $\mathcal{N}(0,4)$ in $[-2,2]$ with $10\%$  nonzero elements. The noise $\omega$ is generated by the standard Gaussian distribution and normalized to the norm $\|\omega\|=10^{-6}$,
\[ b= A{x}^*+\omega.\]
In the numerical examples, we solve the subproblems by the  fast iterative shrinkage-thresholding algorithm  (FISTA) \cite{Beck2009} with the stopping condition
\[ \frac{\|z_k-z_{k-1}\|^2}{\max\{\|z_{k-1}\|,1\}}\leq subtol \]
or the number of iterations exceeds $100$, where $z_k$ is the iterative sequence of FISTA to solve the subproblem and $subtol$ is the precision.
Denote the relative error $Rel = \frac{\|x-x^*\|}{\|{x}^*\|}$ and the residual error $Res = \|Ax-b\|$.
\begin{figure}[htbp]
\centering
\subfigure[$\emph{subtol}=10^{-6}$]{
\begin{minipage}[t]{0.95\linewidth}
\centering
\includegraphics[width=1.55in]{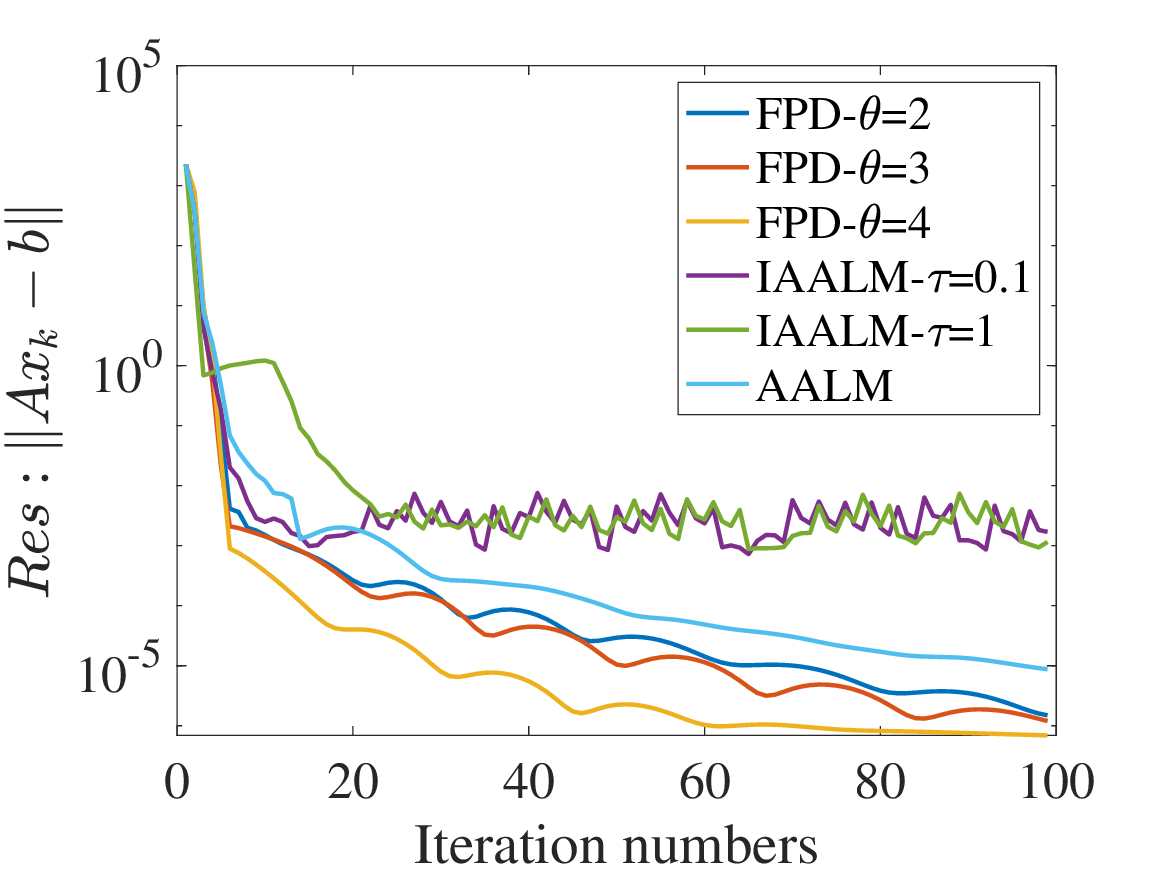}
\centering
\includegraphics[width=1.55in]{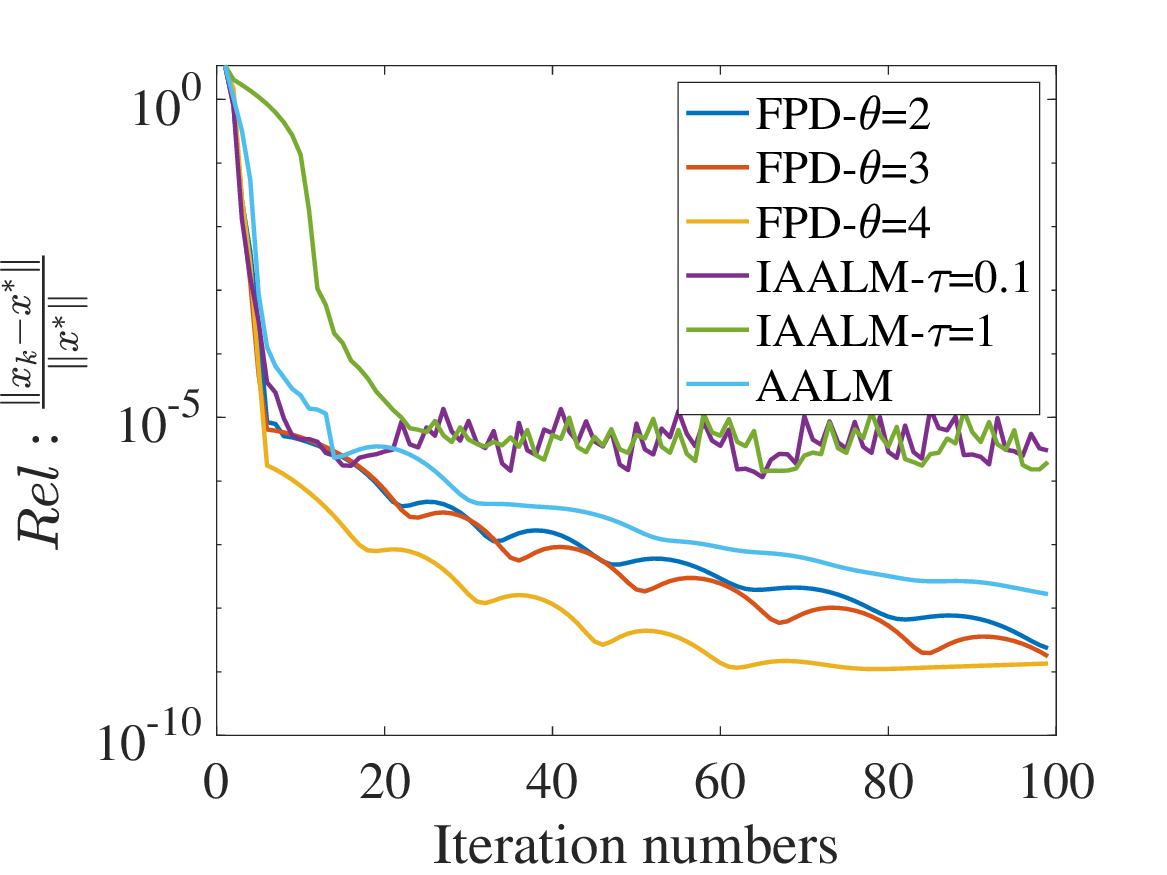}
\end{minipage}%
}%
\\
\subfigure[$\emph{subtol}=10^{-8}$]{
\begin{minipage}[t]{0.95\linewidth}
\centering
\includegraphics[width=1.55in]{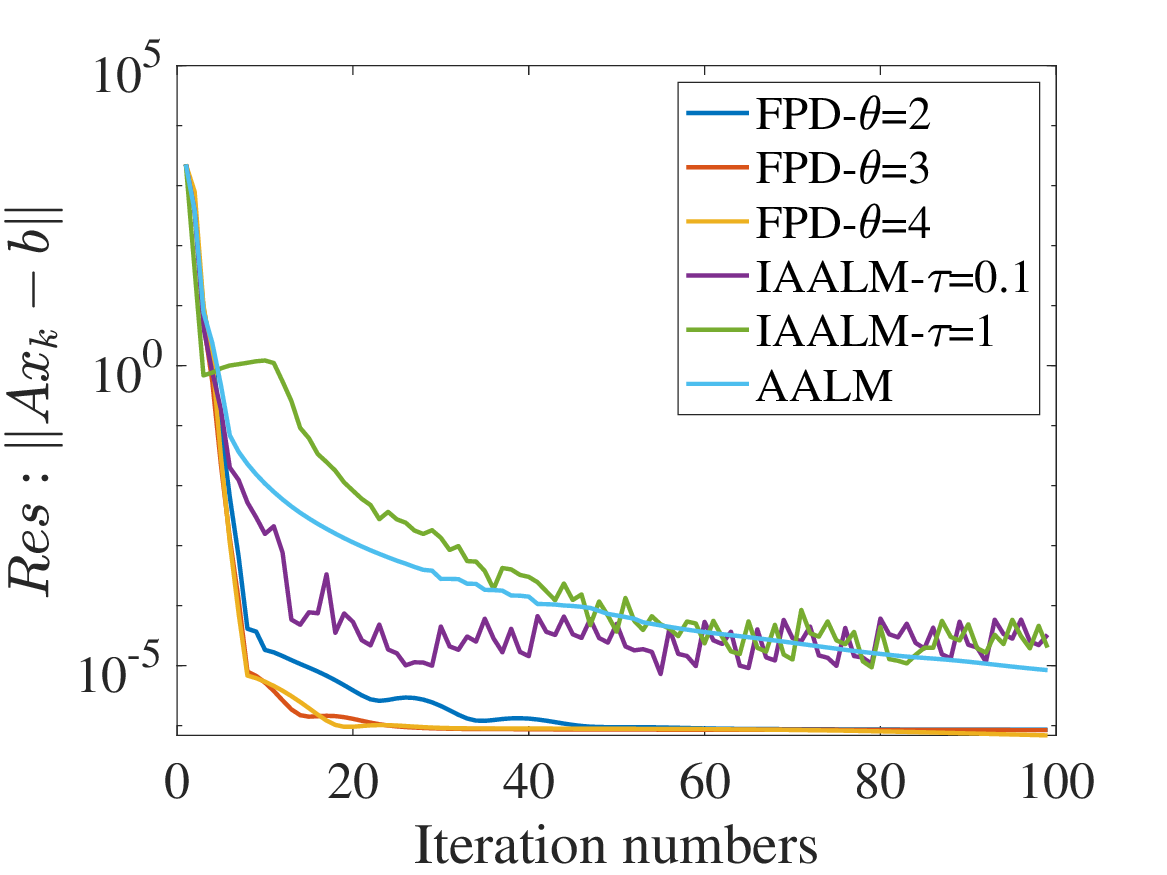}
\centering
\includegraphics[width=1.55in]{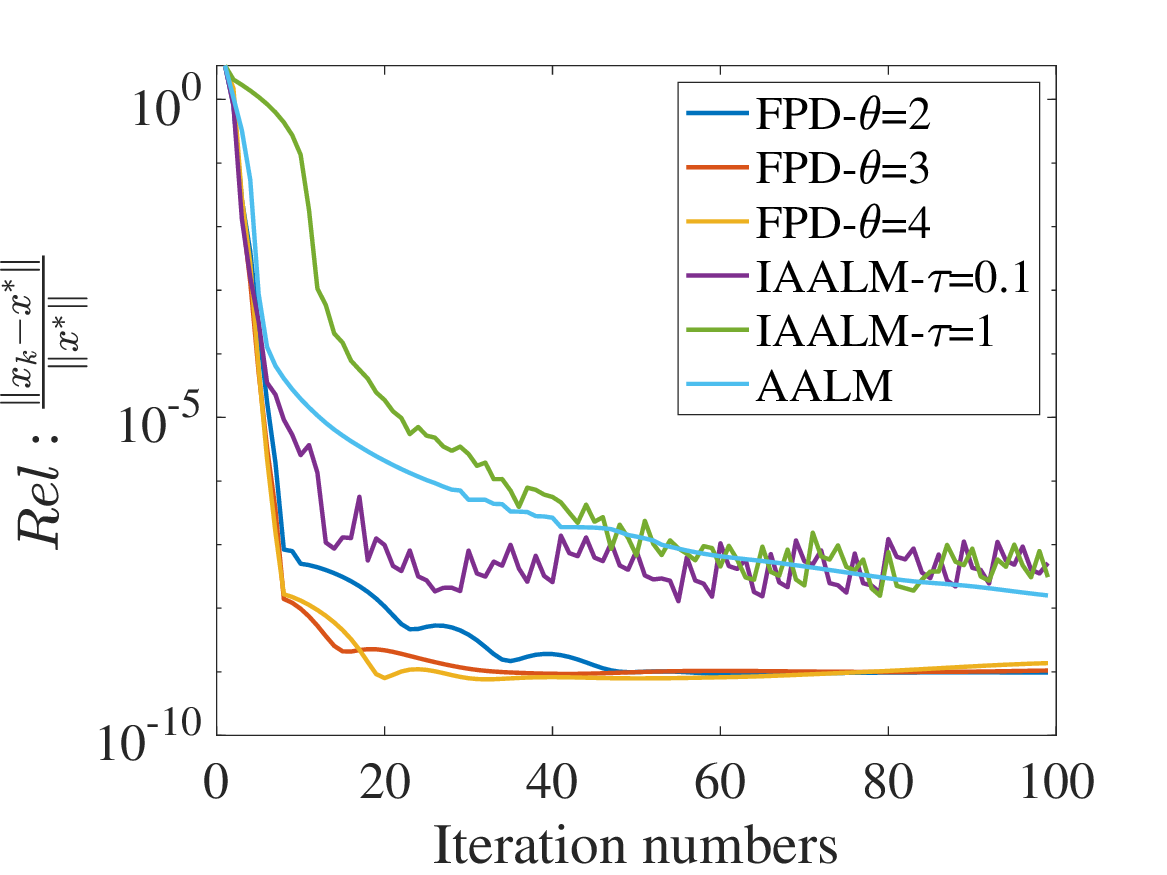}
\end{minipage}%
}%
\\
\subfigure[$\emph{subtol}=10^{-10}$]{
\begin{minipage}[t]{0.95\linewidth}
\centering
\includegraphics[width=1.55in]{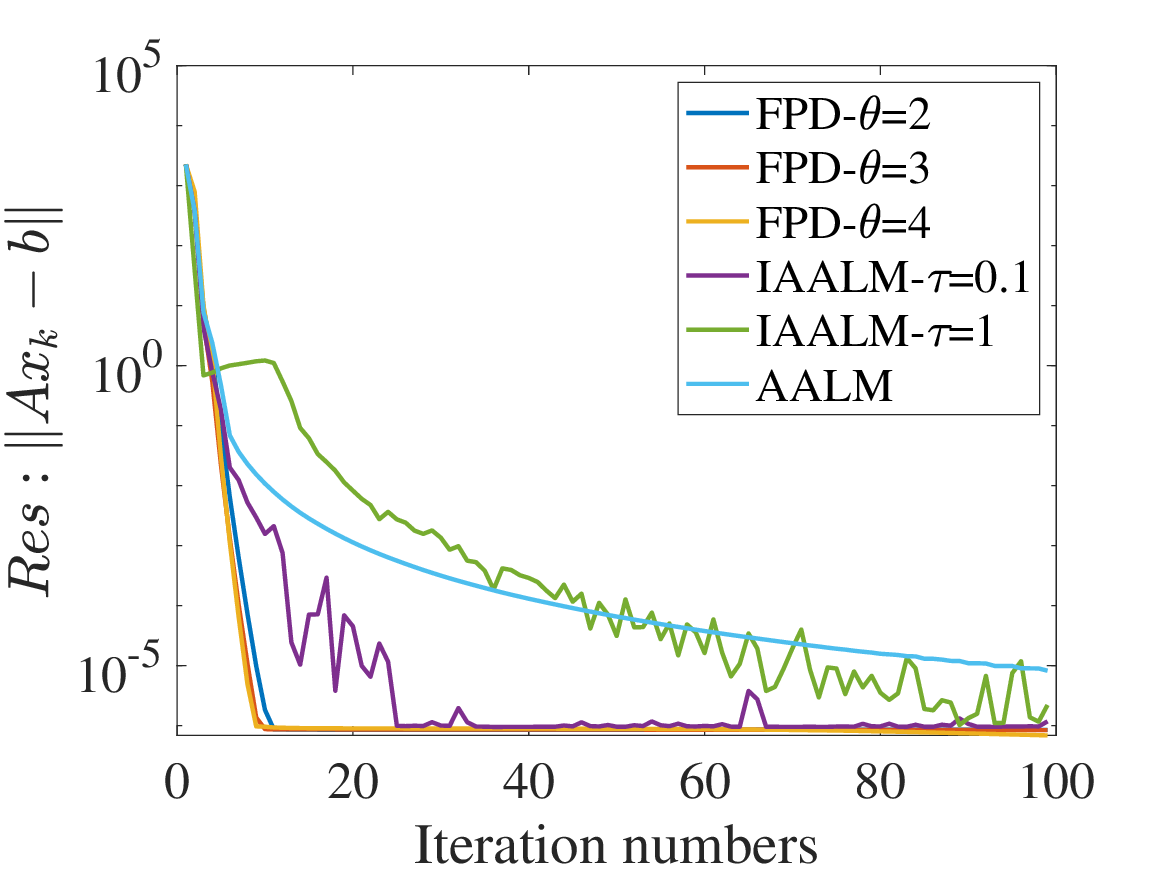}
\centering
\includegraphics[width=1.55in]{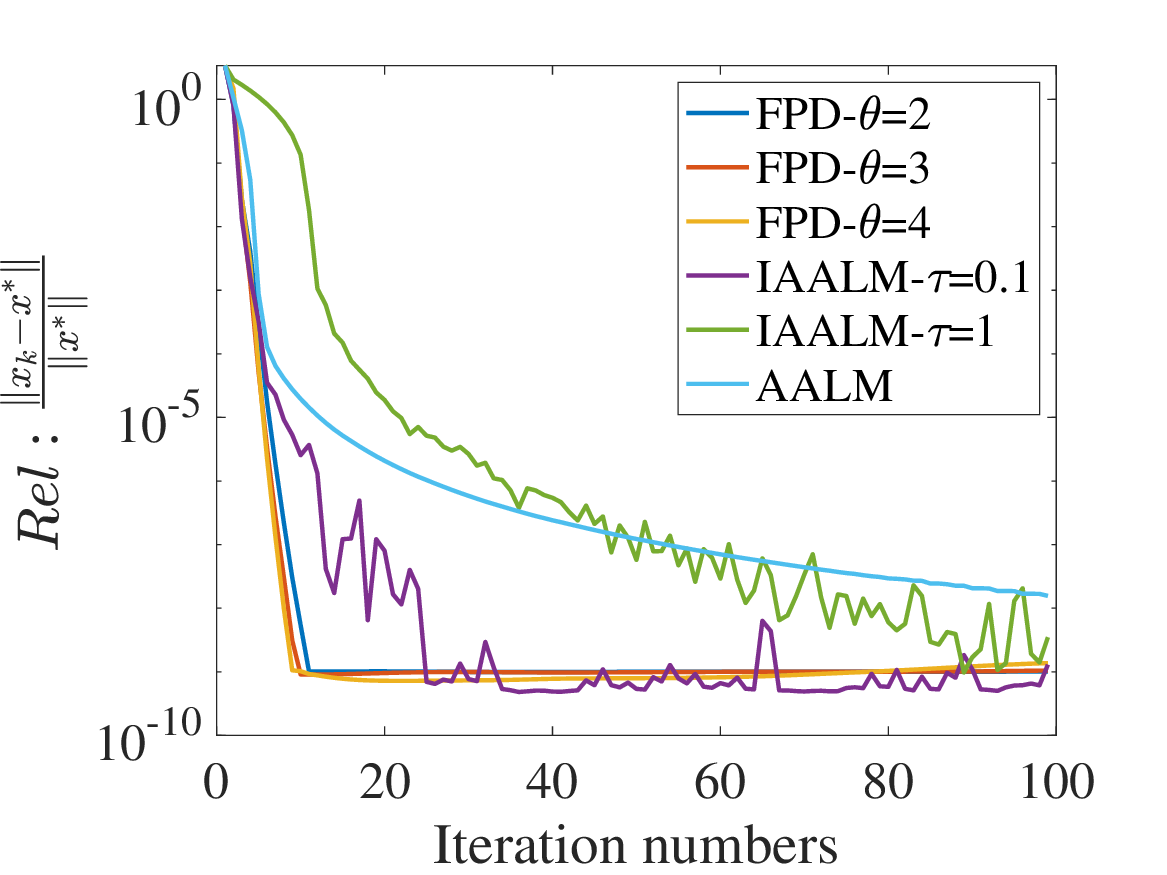}
\end{minipage}%
}%
\caption{Numerical results of FPD, IAALM and AALM under different tolerance of subproblems}\label{fig1}

\end{figure}

 We compare Algorithm \ref{al:al1} (FPD) (update the primal variable by \eqref{subproblem_M}) with IAALM \cite[Algorithm 1]{Kang2015} ($\mathcal{O}(1/k^2)$ convergence rate of the Lagrange residual), and AALM \cite[Algorithm 1]{Xu2017} with adaptive parameters ($\mathcal{O}(1/k^2)$ convergence rate of the objective residual and the feasibility violation). Set the parameters as follows: FPD: $\alpha=50$, $M=\frac{1}{n}Id$, $\beta_0=\frac{0.2}{\theta}$, $\epsilon_k=0$ and 
  \begin{equation*}
   	\beta_{k+1}=\begin{cases}
		\beta_k,  \quad k < \theta-1;\\
		 \frac{k}{k+2-\theta}\beta_k, \quad k\geq \theta-1,
	\end{cases}
\end{equation*} 
where $\theta=\{2,3,4\}$ (In this case, from Lemma \ref{le_new}, the decay rate is $\mathcal{O}(1/k^\theta)$); IAALM: $\tau=\{0.1, 1\}$; AALM: $\gamma=0.1,\ \alpha_k=\frac{2}{k+1},\ \beta_k=\gamma_k=k\gamma,\ P_k=\frac{1}{k}Id$.
 
In Fig. \ref{fig1}, we present the numerical results for various $subtol$ for  first $100$ iterations, which demonstrate the superior performance of Algorithm \ref{al:al1}  over IAALM and AALM under different $subtol$. We can also observe that the larger  $\theta$ is, the better   Algorithm \ref{al:al1} performs.

\section{Conclusion}
By time discretization of the  primal-dual dynamical system  \eqref{dy_pri_dual_a}, we propose an accelerated primal-dual algorithm for the linear equality constrained optimization problem, and prove that the algorithm enjoys the fast convergence rate  $|f(x_k)-f(x^*)|=\mathcal{O}({1}/{k^2\beta_k})$ and $\|Ax_k-b\|=\mathcal{O}({1}/{k^2\beta_k})$.  Further, we prove that the proposed dynamical system  owns a fast convergence properties matching to that of the  algorithm. We exhibit that the known rates from the literature can be obtained  for the second-order dynamical system and the accelerated primal-dual algorithm where only inertial constructions in the sense of Nesterov are needed only for the primal variable.  The numerical experiments demonstrate the validity of acceleration and superior performance of the proposed algorithm over some existing ones.

\appendix
\section{Some auxiliary results}
 The following lemmas have been used in the convergence analysis of the numerical algorithm and the dynamical system.
 
 \begin{lemma}\cite[Lemma 5.14]{Attouch2018MP} \label{le_disc_per}
	Let $\{a_k\}_{k\geq 1}$ and $\{b_k\}_{k\geq 1}$ be two nonnegative sequences. Assume $\sum^{+\infty}_{k=1} b_k< +\infty$ and
	\[ a_k^2 \leq c^2 +\sum^{k}_{j=1}b_j a_j, \quad\forall  k\in\mathbb{N}, \]
 where $c\geq 0$. Then,
	\[\sup_{k\geq 1}a_k\leq c+\sum^{+\infty}_{j=1} b_j<+\infty.\]
	\end{lemma}
	
\begin{lemma}\cite[Lemma A.5]{Brezis1973}\label{le_dy_per}
	Let $\nu:[t_0,T]\to [0,+\infty)$ be integrable and $M\geq 0$. Suppose that $\mu:[t_0,T]\to \mathbb{R}$ is continuous and
	\[ \frac{1}{2}\mu(t)^2\leq \frac{1}{2}M^2+\int^{t}_{t_0}\nu(s)\mu(s)ds\]
	for all $t\in[t_0,T]$. Then, $|\mu(t)|\leq M+\int^t_{t_0}\nu(s)ds$ for all $t\in[t_0,T]$.
\end{lemma}

\begin{lemma}\label{le_A2}
	Let $\{g_k\}_{k\geq k_0}$ be a sequence of vectors in $\mathbb{R}^{n}$ and $\{a_k\}_{k\geq k_0}$ be a sequence in $[0,1)$, where  $k_0\geq 1$. Assume 
	\[\left\|g_{k+1}+\sum_{j=k_0}^ka_j g_j\right\|\leq C, \quad\forall k\geq k_0.\]
	 Then, 
\[ \sup_{k\geq k_0}\|g_k\|<+\infty.\]
\end{lemma}

\begin{pf}
	Define $\{G_k\}_{k\geq k_0}$ be a sequence  of vectors in $\mathbb{R}^{n}$ as 
\begin{equation}\label{eq_defin_G}
	G_k = \rho_k \sum_{j=k_0}^k a_jg_j
\end{equation}
with  $\rho_{k_0} = 1$ and 
\[\rho_{k+1} = \frac{\rho_k}{1-a_{k+1}}, \quad \forall k\geq k_0.\] 
Since $a_k\in[0,1)$,  $\rho_{k+1}-\rho_k= \frac{a_{k+1}\rho_k}{1-a_{k+1}}= \rho_{k+1}a_{k+1}\geq 0$. A direct computation leads to 
\begin{eqnarray*}
	G_{k+1}-G_k &=& \rho_{k+1}\sum_{j=k_0}^{k+1} a_jg_j - \rho_k \sum_{j=k_0}^k a_jg_j\\
	&=& \rho_{k+1}a_{k+1}g_{k+1} + (\rho_{k+1}-\rho_k)\sum_{j=k_0}^{k} a_jg_j\\
	&=& (\rho_{k+1}-\rho_k)(g_{k+1} + \sum_{j=k_0}^{k} a_jg_j),
\end{eqnarray*}
which together with assumption yields
\[\|G_{k+1}-G_k\| \leq C(\rho_{k+1}-\rho_k),\quad \forall k\geq k_0.\]
Using triangle inequality, we get
\begin{eqnarray*}
	\|G_k\| &=&\|G_{k_0}+\sum_{j=k_0}^{k-1}(G_{j+1}-G_j)\| \\
	&\leq& \|G_{k_0}\| +\sum_{j=k_0}^{k-1}\|G_{j+1}-G_j\|\\
	&\leq& \|G_{k_0}\| +C\sum_{j=k_0}^{k-1}(\rho_{j+1}-\rho_j)\\
	&=& \|G_{k_0}\| +C(\rho_k-1)
\end{eqnarray*}
for all $k\geq k_0$. Combining \eqref{eq_defin_G} and $\rho_k\geq 1$ we have
\[\left\|\sum_{j=k_0}^{k} a_j g_j\right\|\leq  \|G_{k_0}\| +C,\quad\forall k\geq k_0.\]
This together with assumption and triangle inequality  implies
\[\sup_{k\geq k_0}\|g_k\|\leq \max\{g_{k_0}, \|G_{k_0}\| +2C\}\leq \|g_{k_0}\|+2C<+\infty.\]
\end{pf}

\begin{lemma}\label{le_new}
	Let $\{\gamma_k\}_{k\geq 1}\subset (0,+\infty)$ be a positive sequence such that  $\gamma_{k+1} = (1+\frac{a}{k+b})\gamma_{k}$  for any $k\geq 1$, where $a>0$ and $b\geq 0$.  Then, there exist $\mu_1>0$ and $\mu_2>0$ such that 
	\[ \mu_1 k^a\leq  \gamma_k\leq \mu_2 k^a.\]
\end{lemma}
\begin{pf}
	Define $\varphi:[1,+\infty)\to(0,+\infty)$ as
	\[\varphi(t) = (k+1-t)\gamma_k+(t-k)\gamma_{k+1}=\frac{k+b+a(t-k)}{k+b}\gamma_k\]
for any $k\leq t<k+1$.
	It is easy to verify that $\varphi(t)$ is a positive, piecewise linear and nondecreasing function. By the definition of $\gamma_k$, we can compute
\[\dot{\varphi}(t) = \gamma_{k+1}-\gamma_k = \frac{a}{k+b}\gamma_k =\frac{a}{k+b+a(t-k)}\varphi(t) \]
for any $t\in(k,k+1)$. It yields
\[\frac{\dot{\varphi}(t)}{\varphi(t)} = \frac{a}{(1-a)k+b+at},
\quad\forall t\in(k,k+1).\]
Then for any $t\geq 1$, we have 
\begin{equation*}
	\begin{cases}
		\frac{a}{t+b}\leq\frac{\dot{\varphi}(t)}{\varphi(t)} \leq \frac{a}{t+a+b-1}, & 0\leq a \leq 1  ,\\
		 \frac{a}{t+a+b-1}\leq \frac{\dot{\varphi}(t)}{\varphi(t)} \leq \frac{a}{t+b}, & a>1.
	\end{cases}
\end{equation*}
Since $\varphi(t)$ is a piecewise linear function,  integrating the above inequalities over $[0,t]$, we have 
\[a\ln(t+b)+C_1\leq \ln \varphi(t) \leq a\ln(t+a+b-1)+C_2\]
as $0\leq a \leq 1$,
and 
\[a\ln(t+a+b-1)+C_1\leq \ln \varphi(t)\leq a\ln(t+b)+C_2\]
as $a>1$, where $C_1$ and $C_2$ are two constant.
It follows that
\begin{equation*}
	\begin{cases}
		e^{C_1}(t+b)^a\leq  \varphi(t) \leq e^{C_2}(t+a+b-1)^a,\ & 0\leq a \leq 1,\\
		 e^{C_1}(t+a+b-1)^a\leq  \varphi(t)\leq e^{C_2}(t+b)^a,\  & a>1.
	\end{cases}
\end{equation*}
As a result, for any $a>0$ there exists $\mu_1>0$ and $\mu_2>0$ such that $\mu_1 t^a\leq  \varphi(t)\leq \mu_2 t^a$. Letting $t=k$ leads to the result.
\end{pf}

\begin{lemma}\label{le_A4}
	Assume that $g:[t_0,+\infty)\to\mathbb{R}^n$ is a continuous  function, $a:[t_0,+\infty)\to[0,+\infty)$ is a continuous  function,  $t_0>0$, and $C\geq 0$. If

\begin{equation}\label{lm7}
\left\|g(t)+\int^t_{t_0}a(s){g}(s)ds\right\|\leq C, \qquad \forall t\geq t_0,
\end{equation}
	then 
	\[\sup_{t\geq t_0}\|{g}(t)\|<+\infty.  \]
\end{lemma}

\begin{pf}
	Define $G:[t_0,+\infty)\to \mathbb{R}^n$ by
	\begin{equation}\label{def_G}
		G(t)=e^{\int^t_{t_0}a(s)ds}\cdot\int^t_{t_0}a(s){g}(s)ds.
	\end{equation}
From \eqref{lm7} and \eqref{def_G}, we have	
\[\left\|\dot{G}(t)\right\|\leq Ca(t)e^{\int^t_{t_0}a(s)ds}.\]
Note that $G(t_0)=0$. It follows that 
\begin{eqnarray*}
	\|G(t)\| &=& \left\|\int^t_{t_0}\dot{G}(s)ds\right\|\leq \int^t_{t_0}\left\|\dot{G}(s)\right\|ds\\
	&\leq&   \int^t_{t_0}Ca(s)e^{\int^s_{t_0}a(h)dh}ds = Ce^{\int^t_{t_0}a(s)ds}-C.
\end{eqnarray*}
This together with \eqref{def_G} yields 
\[\left\|\int^t_{t_0}a(s){g}(s)ds\right\| \leq C,\quad \forall t\geq t_0.\]
From  \eqref{lm7} and triangle inequality, we get
\[\sup_{t\geq t_0}\|{g}(t)\|\leq  2C<+\infty.\]
\end{pf}

\begin{ack}                               
The authors would like to thank the reviewers and the editors for their helpful comments and suggestions,  which  significantly improve  the quality of this paper.  In particular, we would like to thank one reviewer for sharing Lemma 4 and Lemma 5 with us, which have made the results better.
\end{ack}
%

\end{document}